\documentclass[11pt]{article}
\usepackage{exscale,relsize}
\usepackage{amsmath}
\usepackage{amsfonts}
\usepackage{amssymb}
\usepackage{calc}
\usepackage{theorem}
\usepackage{pifont}      
\newtheorem{theorem}{Theorem}[section]

\newtheorem{fact}[theorem]{Fact}
\newtheorem{corollary}[theorem]{Corollary}

\newtheorem{definition}[theorem]{Definition}

\theoremstyle{plain}{\theorembodyfont{\rmfamily}
}
\theoremstyle{plain}{\theorembodyfont{\rmfamily}
}
\theoremstyle{plain}{\theorembodyfont{\rmfamily}
}
\theoremstyle{plain}{\theorembodyfont{\rmfamily}
\newtheorem{example}[theorem]{Example}}
\theoremstyle{plain}{\theorembodyfont{\rmfamily}
\newtheorem{remark}[theorem]{Remark}}
\theoremstyle{plain}{\theorembodyfont{\rmfamily}
}

\newcommand{\gra}{\ensuremath{\operatorname{gra}}}
\newcommand{\scal}[2]{\langle{{#1},{#2}}\rangle}
\newcommand{\Id}{\ensuremath{\operatorname{Id}}}
\newcommand{\To}{\ensuremath{\rightrightarrows}}
\newcommand{\menge}[2]{\big\{{#1} \mid {#2}\big\}}
\newcommand{\ran}{\ensuremath{\operatorname{ran}}}
\newcommand{\RR}{\ensuremath{\mathbb R}}

\newcommand{\inte}{\ensuremath{\operatorname{int}}}

\newcommand{\clconv}{\ensuremath{\overline{\operatorname{conv}}\,}}
\newcommand{\cldom}{\ensuremath{\overline{\operatorname{dom}}\,}}
\newcommand{\clran}{\ensuremath{\overline{\operatorname{ran}}\,}}
\newcommand{\conv}{\ensuremath{\operatorname{conv}}}
\newcommand{\dom}{\ensuremath{\operatorname{dom}}}
\newcommand{\reli}{\ensuremath{\operatorname{ri}}}
\newcommand{\aff}{\ensuremath{\operatorname{aff}}}
\newcommand{\ospan}{\ensuremath{\operatorname{span}}}

\begin{document}

\title{Firmly nonexpansive and Kirszbraun-Valentine extensions:\\
a constructive approach via monotone operator theory}

\author{
Heinz H.\ Bauschke\thanks{Mathematics, Irving K.\ Barber School,
UBC Okanagan, Kelowna, British Columbia V1V 1V7, Canada. E-mail:
\texttt{heinz.bauschke@ubc.ca}.}~~ and ~~ Xianfu
Wang\thanks{Mathematics, Irving K.\ Barber School, UBC Okanagan,
Kelowna, British Columbia V1V 1V7, Canada. E-mail:
\texttt{shawn.wang@ubc.ca}.} }

\date{July 8, 2008\\
~~\\\emph{Dedicated to Alex Ioffe and Simeon Reich,\\ on the occasion of
their 70th and 60th birthdays}}

\maketitle

\begin{abstract}
\noindent
Utilizing our recent proximal-average based
results on the constructive extension of monotone operators,
we provide a novel approach to the celebrated Kirszbraun-Valentine Theorem
and to the extension of firmly nonexpansive mappings.
\end{abstract}



\section{Introduction}

Throughout, $X$ is a real Hilbert space with inner product
$\scal{\cdot}{\cdot}$ and associated norm $\|\cdot\|$.

\begin{definition}
Let $S \subseteq X$ and let $T\colon S\to X$.
Then $T$ is \emph{nonexpansive} if
\begin{equation}
(\forall x\in S)(\forall y\in S)\quad
\|Tx-Ty\|\leq\|x-y\|.
\end{equation}
\end{definition}
See, e.g.\ \cite{GK90,GR84} for further information on nonexpansive
mappings.
Let us recall the celebrated Kirszbraun-Valentine Theorem
(see \cite{K34,V43,V45}), which
states that every nonexpansive mapping can be extended to the
entire space.

\begin{fact}[Kirszbraun-Valentine]
Let $S\subseteq X$ and let $T\colon S\to X$ be nonexpansive.
Then there exists a mapping $\widetilde{T}\colon X\to X$ such that
$\widetilde{T}$ is nonexpansive and $\widetilde{T}|_S = T$.
\end{fact}

In this note, we present constructive approaches to the Kirszbraun-Valentine
extension theorem for nonexpansive mappings. Our main tool is the proximal average,
which is used to provide an explicit
maximal monotone extension of any monotone operator.
In Section~\ref{monoextension}, we give maximal monotone extensions
with ``small domains'' under some mild constraint qualifications.
Our extension results for firmly nonexpansive and for
nonexpansive mappings are presented in Section~\ref{main}.
In some cases, it is possible to provide constructive extensions
with an optimally localized range.

\section{Firmly Nonexpansive Mappings and Monotone Operators}\label{monoextension}

\begin{definition}
Let $S\subseteq X$ and let $F\colon S\to X$.
Then $F$ is \emph{firmly nonexpansive} if
\begin{equation}
(\forall x\in S)(\forall y\in S)\quad
\|Fx-Fy\|^2\leq\scal{Fx-Fy}{x-y}.
\end{equation}
\end{definition}

The Cauchy-Schwarz inequality implies
that every firmly nonexpansive operator is nonexpansive.
The next result, which is well known
(see, e.g., \cite[Theorem~12.1]{GK90}),
provides a bijection between nonexpansive and firmly nonexpansive
mappings.

\begin{fact} \label{f:nonexpvsfirm}
Let $S\subseteq X$, let $T\colon S\to X$, and let $F\colon S\to X$.
Suppose that $F = \tfrac{1}{2}\Id+\tfrac{1}{2}T$.
Then $T$ is nonexpansive $\Leftrightarrow$ $F$ is firmly nonexpansive.
\end{fact}

\begin{remark} \label{r:last1}
Note that (see also \cite{LS87})
the linear relationship between $\gra T$
and $\gra F$. In fact,
the linear operator 
$(x,y)\mapsto (x,\tfrac{1}{2}x+\tfrac{1}{2}y)$
provides a bijection from $\gra T$ to $\gra F$
with inverse $(x,y)\mapsto (x,2y-x)$.
\end{remark}

\begin{definition}
Let $A\colon X\To X$, i.e., $A$ is a set-valued
operator from $X$ to the power set of $X$.
Denote the \emph{graph} of $A$ by $\gra A :=
\menge{(x,x^*)\in X\times X}{x^*\in Ax}$.
Then $A$ is \emph{monotone} if
\begin{equation}
\big(\forall (x,x^*)\in\gra A\big)\big(\forall (y,y^*)\in \gra A\big)\quad
\scal{x-y}{x^*-y^*} \geq 0.
\end{equation}
If $(x,x^*)\in X\times X$ and the operator with graph
$\{(x,x^*)\}\cup \gra A$ is monotone, then $(x,x^*)$ is
\emph{monotonically related} to $\gra A$.
If $A$ is monotone and every proper extension of $A$ fails to be
monotone, then $A$ is \emph{maximal monotone}.
The \emph{inverse} operator $A^{-1}$ is defined via
$\gra A^{-1} := \menge{(x^*,x)\in X\times X}{x^*\in Ax}$.
\end{definition}

Maximal monotone operators play a critical role in modern
Analysis and Optimization; see, e.g.,
\cite{BI07,RW98,S98,S08,Z02}.

The following result,
brought out fully by Eckstein and Bertsekas
\cite{EB92}, has its roots in the seminal works by Minty \cite{M62}
and by Rockafellar \cite{R76}.

\begin{fact} \label{f:EB}
Let $A\colon X\To X$, let $S\subseteq X$, and let
$F\colon S\to X$.
Suppose that $F = (A+\Id)^{-1}$; equivalently, that $A = F^{-1}-\Id$.
Then the following hold.
\begin{enumerate}
\item $A$ is monotone $\Leftrightarrow$ $F$ is firmly nonexpansive.
\item $A$ is maximal monotone $\Leftrightarrow$
$F$ is firmly nonexpansive and $S=X$.
\end{enumerate}
\end{fact}

\begin{remark} \label{r:last2}
Using notation of Fact~\ref{f:EB}, we recall
the linear relationship (see \cite{LS87})
between the graphs of $F$ and $A$.
Indeed, the linear operator
$(x,y)\mapsto (y,x-y)$
provides a bijection from $\gra F$ to $\gra A$
with inverse $(x,y)\mapsto (x+y,x)$.
Recalling Remark~\ref{r:last1}, we observe the linear
relationship from $\gra T$ to $\gra A$ via
$(x,y)\mapsto \tfrac{1}{2}(x+y,x-y)$ with
inverse $(x,y)\mapsto (x+y,x-y)$.
\end{remark}

\begin{corollary} \label{c:neato}
Let $F\colon X\to X$ be firmly nonexpansive.
Then $\clran F$ is convex.
\end{corollary}
\begin{proof}
Set $A := T^{-1}-\Id$ and observe that $\dom A = \ran T$.
The conclusion now follows from Fact~\ref{f:EB}.(ii)
and \cite[Theorem~18.6]{S98}.
\end{proof}

\begin{remark}
Suppose that $X=\RR^2$,
let $P_A$ the projector onto the line $\RR\times\{1\}$
and let $P_B$ be the projector onto the closed unit ball.
Then $P_A$ and $P_B$ are both (firmly) nonexpansive
(see, e.g., \cite[Theorem~12.2]{GK90}).
Set $T := P_BP_A$.
Then $T$ is a nonexpansive mapping defined on the
entire Euclidean plane; however, $\clran T$ is
equal to the the closed upper half circle, which is not convex.
\end{remark}

The following result, originally obtained with
the help of the \emph{proximal average},
the \emph{Fitzpatrick function}, and
other tools from \emph{Convex Analysis}, is of key importance.
We refer the reader to
\cite{BGLW08,BLT08,BLW07,BMR04,BW07b}
for further information on the proximal average,
to \cite{BS02,F88,ML01,P04,S08} and the references therein for
results on the Fitzpatrick function,
and to
\cite{R76,RW98,Z02} for
the basic theory of Convex Analysis.
The notation $\widetilde{A}$ for the maximal monotone extension
of a monotone operator $A\colon X\To X$ will be
used from now on.

\begin{fact} (See \cite[Fact~5.6 and Theorem~5.7]{BW07b}.)
\label{f:BW}
Let $A\colon X\To X$ be monotone.
Recall that the \emph{Fitzpatrick function}
of $A$ is the function on $X\times X$ defined by
\begin{equation} \label{e:defofF_A}
\Phi_A\colon (x,x^*)\mapsto \sup_{(a,a^*)\in\gra A}\big(
\scal{x}{a^*} + \scal{a}{x^*} - \scal{a}{a^*}\big),
\end{equation}
with \emph{Fenchel conjugate}
\begin{equation}
\Phi_A^*\colon (y^*,y)\mapsto \sup_{(x,x^*)\in X\times X}
\big(\scal{x}{y^*} + \scal{y}{x^*}-\Phi_A(x,x^*)\big).
\end{equation}
Set
\begin{equation} \label{e:defofPsi}
\Psi_A\colon (x,x^*) \mapsto
\min_{(x,x^*)=\tfrac{1}{2}(y+z,y^*+z^*)}\left(
\tfrac{1}{2}\Phi_A(y,y^*) + \tfrac{1}{2}\Phi_A^*(z^*,z) + \tfrac{1}{8}
\big(\|y-z\|^2+\|y^*-z^*\|^2\big)\right),
\end{equation}
which is the \emph{proximal average} between $\Phi_A$
and the (transpose of the) $\Phi_A^*$,
and define $\widetilde{A}\colon X\To X$ via
\begin{equation} \label{e:defofAtilde}
\gra \widetilde{A} = \menge{(x,x^*)\in X\times X}{(x^*,x)\in\partial
\Psi_A(x,x^*)},
\end{equation}
where ``$\partial$''
denotes the \emph{subdifferential operator} from Convex Analysis.
Then the following hold.\hfill
\begin{enumerate}
\item \label{f:BWi}
$(\forall (x,x^*)\in X\times X)$
$\Psi_A(x,x^*) \geq \scal{x}{x^*}$.
\item \label{f:BWii}
$(\forall (x,x^*)\in X\times X)$
$x^*\in\widetilde{A}x \Leftrightarrow
\Psi_A(x,x^*) = \scal{x}{x^*}$.
\item \label{f:BWiii}
$\widetilde{A}$ is a maximal monotone extension of $A$.
\end{enumerate}
\end{fact}


It is convenient (see, e.g., Example~\ref{ex:wednesday} below)
to be able to use
the following alternative description of $\widetilde{A}$.

\begin{theorem} \label{t:wednesday}
Let $A:X\To X^*$ be monotone. Define $B:X\To X^*$ via
\begin{equation}
x^*\in Bx \quad \Leftrightarrow \quad
(x,x^*) =\tfrac{1}{2}(x_1+x_2,x_1^*+x_2^*),
\end{equation}
where
\begin{equation} \label{e:neun}
(x_{2}^*,x_{2})\in\partial
\Phi_A(x_{1},x_{1}^*)\quad\text{and}
\quad  x_1-x_{2}=x_{1}^*-x_{2}^*.
\end{equation}
Then
$B = \widetilde{A}$.
\end{theorem}
\begin{proof}
Take $(x,x^*)\in X\times X$.

``$ \gra \widetilde{A}\subseteq \gra B$'':
Suppose that $(x,x^*)\in \gra \widetilde{A}$.
Then $\Psi_{A}(x,x^*)=\langle x^*,x\rangle$ by
Fact~\ref{f:BW}.\ref{f:BWii} and,
by \eqref{e:defofAtilde},
there exist $(x_1,x_1^*)$ and $(x_2,x_2^*)$ in
$X\times X$ such that
\begin{equation} \label{e:wednesday:b}
(x,x^*)=\tfrac{1}{2}(x_{1}+x_{2},x_{1}^*+x_{2}^*)
\end{equation}
and
\begin{equation}\label{e:achieved}
\tfrac{1}{4}\scal{x_{1}+x_{2}}{x_{1}^*+x_{2}^*}  =
\tfrac{1}{2} \Phi_A(x_{1},x_{1}^*)+\tfrac{1}{2}\Phi_A^*(x_{2}^*,x_{2})
 +\tfrac{1}{8}\big(\|x_1-x_{2}\|^2+\|x_{1}^*-x_{2}^*\|^2\big).
\end{equation}
The Fenchel-Young inequality yields
$\Phi_A(x_{1},x_{1}^*)+\Phi_A^*(x_{2}^*,x_{2})\geq
\scal{(x_1,x_1^*)}{(x_2^*,x_2)} = \langle
x_{1},x_{2}^*\rangle+\langle x_{1}^*,x_{2}\rangle$.
Hence, \eqref{e:achieved} results in
\begin{equation}
\tfrac{1}{4}\scal{x_{1}+x_{2}}{x_{1}^*+x_{2}^*}  \geq
\tfrac{1}{2}\langle x_{1},x_{2}^*\rangle +\tfrac{1}{2}\langle
x_1^*,x_{2}\rangle
 +\tfrac{1}{8}\big(\|x_1-x_{2}\|^2+\|x_{1}^*-x_{2}^*\|^2\big);
\end{equation}
equivalently,
\begin{equation}
\tfrac{1}{2}\|x_1-x_2\|^2 + \tfrac{1}{2}\|x_1^*-x_2^*\|^2
\leq \scal{x_1-x_2}{x_1^*-x_2^*}.
\end{equation}
Thus, $x_1^*-x_2^* \in \partial \big(\tfrac{1}{2}\|\cdot\|^2\big)(x_1-x_2)$
and hence
\begin{equation} \label{e:wednesday:a}
x_1^*-x_2^*=x_1-x_2.
\end{equation}
In view of
\eqref{e:achieved} and \eqref{e:wednesday:a}, it follows that
\begin{equation}
\tfrac{1}{2}\scal{x_{1}+x_{2}}{x_{1}^*+x_{2}^*}  =
\Phi_A(x_{1},x_{1}^*)+\Phi_A^*(x_{2}^*,x_{2})
 +\tfrac{1}{2}\scal{x_1-x_2}{x_1^*-x_2^*},
\end{equation}
i.e.,
$\Phi_A(x_1,x_1^*) + \Phi_A^*(x_2^*,x_2)
= \scal{x_1}{x_2^*} + \scal{x_2}{x_1^*}
= \scal{(x_1,x_1^*)}{(x_2^*,x_2)}$
and hence
\begin{equation}
\label{e:wednesday:c}
(x_{2}^*,x_{2})\in\partial \Phi_A(x_{1},x_{1}^*).
\end{equation}
Combining \eqref{e:wednesday:b}, \eqref{e:wednesday:c},
and \eqref{e:wednesday:a}, we see that $(x,x^*)\in\gra B$.

``$\gra B \subseteq \gra \widetilde{A}$'':
 Suppose that $(x,x^*)\in\gra B$.
 Then there exist $(x_1,x_1^*)$ and $(x_2,x_2^*)$ in $X\times X$.
 such that \eqref{e:wednesday:b} and \eqref{e:neun} hold.
 Note that \eqref{e:neun} yields
\begin{equation} \label{e:wednesday:d}
\Phi_A(x_{1},x_{1}^*)+\Phi_A^*(x_{2}^*,x_{2})=\langle
 x_{1},x_{2}^*\rangle +\langle x_{1}^*,x_{2}\rangle
\end{equation}
 and
\begin{equation} \label{e:wednesday:e}
\tfrac{1}{2}\|x_{1}-x_{2}\|^2 +
\tfrac{1}{2}\|x_{1}^*-x_{2}^*\|^2 =
\langle x_1-x_{2},x_{1}^*-x_{2}^*\rangle.
\end{equation}
Using \eqref{e:wednesday:b}, \eqref{e:wednesday:d},
and \eqref{e:wednesday:e}, we obtain
\begin{align}
\langle x,x^*\rangle &= \tfrac{1}{4} \langle
x_{1}+x_{2},x_{1}^*+x_{2}^*\rangle \label{longset1}\\
&= \tfrac{1}{2}\Phi_A(x_{1},x_{1}^*)+
\tfrac{1}{2}\Phi_A^*(x_{2}^*,x_{2})+
\tfrac{1}{4}\langle x_{1}^*-x_{2}^*,x_{1}-x_{2}\rangle\\
&= \tfrac{1}{2}\Phi_A(x_{1},x_{1}^*)+
\tfrac{1}{2}\Phi_A^*(x_{2}^*,x_{2})+
\tfrac{1}{8}\big(\|x_{1}-x_{2}\|^2 + \|x_{1}^*-x_{2}^*\|^2\big).
\label{longset3}
\end{align}
Using this, \eqref{e:defofPsi}, and
Fact~\ref{f:BW}.\ref{f:BWi}, we estimate
\begin{align}
\scal{x}{x^*} &=
\tfrac{1}{2}\Phi_A(x_{1},x_{1}^*)+
\tfrac{1}{2}\Phi_A^*(x_{2}^*,x_{2})+
\tfrac{1}{8}\big(\|x_{1}-x_{2}\|^2 + \|x_{1}^*-x_{2}^*\|^2\big)\\
&\geq \Psi_A(x,x^*)\notag\\
&\geq \scal{x}{x^*},\notag
\end{align}
and we see that $\Psi_{A}(x,x^*) = \scal{x}{x^*}$.
By Fact~\ref{f:BW}\ref{f:BWii}, $(x,x^*)\in\gra\widetilde{A}$.
\end{proof}

\begin{example}\label{ex:simplest}
Let $(a,a^*)\in X\times X$ and let $A\colon X\To X$
be given by $\gra A = \{(a,a^*)\}$.
Then $\gra \widetilde{A} = \{(x,x+a^*-a)\}_{x\in X}$.
\end{example}
\begin{proof}
In view of \eqref{e:defofF_A}, we see that
\begin{equation}
\Phi_A\colon X\to\RR\colon (x,x^*) \mapsto
\scal{x}{a^*} + \scal{a}{x^*} - \scal{a}{a^*}
=\scal{(x,x^*)}{(a^*,a)} - \scal{a}{a^*}.
\end{equation}
Since $\partial \Phi_A \equiv (a^*,a)$,
Theorem~\ref{t:wednesday} implies
\begin{equation}
\gra \widetilde{A} = \menge{\tfrac{1}{2}(x_1+a,x_1-a+2a^*)}{x_1\in X}.
\end{equation}
The change of variable $x=\tfrac{1}{2}(x_1+a)$
now yields the result.
\end{proof}

Note that in Example~\ref{ex:simplest}
the domain of the extension $\widetilde{A}$ is the entire
space $X$ while the domain of the given operator
is a singleton.
This raises the question on finding extensions with the smallest possible
domain.
Let $A\colon X\To X$ be monotone and set $D := \clconv \dom A$.
Denote the \emph{normal cone} operator to $D$ by $N_D$; i.e.,
$N_D = \partial\iota_D$ so that $\dom N_D=D$ and $(\forall x\in D)$
$N_D(x) = \menge{x^*\in X}{\sup\scal{D-x}{x^*}\leq 0}$.
Since $(\forall x\in D)$ $0\in N_D(x)$, the operator
\begin{equation} \label{e:A+N_D}
A + N_D
\end{equation}
is monotone extension of $A$. However, $A+N_D$ may fail to be
maximal monotone: consider, e.g., the case when
$A$ is the zero operator
restricted to the \emph{open} unit ball.
The following result will aid us in our quest to provide
a sufficient condition for \eqref{e:A+N_D} to be maximal
monotone.

\begin{fact} (See \cite[Theorem~2.14 and Theorem~2.15]{BW07b}.)
\label{f:thursday}
Let $A\colon X\To X$ be monotone, and set $D := \clconv \dom A$.
Then the following hold. \hfill
\begin{enumerate}
\item \label{f:thursdayi}
$(x,x^*)\in X\times X$ is monotonically related to $\gra(A+N_D)$
$\;\Leftrightarrow\;$
\begin{equation}
\text{$(x,x^*)$ is monotonically related to $\gra A$}\quad
\text{and}\quad
x\in \bigcap_{a\in \dom A}\big(a + T_D(a)\big),
\end{equation}
where $T_D(a) = N_D^\ominus(a)$ is the \emph{polar}
(negative dual) cone of $N_D(a)$.
\item \label{f:thursdayii}
If $\conv\dom A$ is closed, then $\bigcap_{a\in\dom A}(a+
T_D(a)) = D$.
\end{enumerate}
\end{fact}

\begin{theorem} \label{t:thursday}
Let $A\colon X\To X$ be monotone, set $D := \clconv\dom A$ and
$A_D := A + N_D$. Then the following hold.
\begin{enumerate}
\item  \label{t:thursdayi}
$\widetilde{A_D}$ is a maximal monotone extension of $A$,
and
\begin{equation} \label{e:thursday:a}
\cldom\widetilde{A_D} = \clconv\dom\widetilde{A_D}
\subseteq \bigcap_{a\in \dom A}
\big(a+T_{D}(a)\big).
\end{equation}
\item \label{t:thursdayii}
If $\conv\dom A$ is closed, then $\cldom \widetilde{A_D} = D$.
\end{enumerate}
\end{theorem}
\begin{proof}
Note that $\dom A_D = \dom A \cap \dom N_D = \dom A \cap D = \dom A$.

\ref{t:thursdayi}:
By Fact~\ref{f:BW}.\ref{f:BWiii},
$\widetilde{A_D}$ is a maximal monotone extension of $A_D$.
On the other hand, $A_D$ is a monotone extension of $A$.
Altogether, $\widetilde{A_D}$ is a maximal monotone
extension of $A$.
Fact~\ref{f:thursday}.\ref{f:thursdayi} yields
$\dom\widetilde{A_{D}}\subseteq\bigcap_{a\in\dom A}(a+T_D(a))$.
Since the last intersection is closed and convex,
we obtain the inclusion in
\eqref{e:thursday:a}.
Finally, the maximal monotonicity of $\widetilde{A_D}$ and
\cite[Theorem~18.6]{S98} imply that
$\clconv\dom\widetilde{A_D} = \cldom\widetilde{A_D}$.

\ref{t:thursdayii}: Combining
\ref{t:thursdayi} and Fact~\ref{f:thursday}.\ref{f:thursdayii},
we deduce that
\begin{equation}
D = \clconv \dom A \subseteq
\clconv\dom \widetilde{A_D} = \cldom\widetilde{A_D}
\subseteq \bigcap_{a\in\dom A}\big(a+T_D(a)\big) = D.
\end{equation}
The proof is complete.
\end{proof}

\begin{remark}
In Theorem~\ref{t:thursday}.\ref{t:thursdayii},
the set $\conv\dom A$ is closed whenever one
of the following holds.
\begin{enumerate}
\item $\dom A$ is a finite subset of $X$.
\item $X$ is finite-dimensional, and $\dom A$ is closed and bounded.
\end{enumerate}
\end{remark}

We conclude this section with a second approach that
tailored to finite-dimensional spaces.

\begin{theorem} \label{t:george}
Suppose that $X$ is finite-dimensional, and let $A\colon X\To X$
be monotone.
Suppose that (after translation if necessary) $0\in \reli \clconv \dom
A=\reli\conv\dom A$, i.e., $0\in \inte_Y \clconv \dom A$,
where $Y:=\aff \dom A=\ospan \dom A$ is a closed subspace of $X$.
Let $P\colon X\to X$ be the linear orthogonal projector onto $Y$,
and let $Q\colon X\to Y\colon x\mapsto Px$.
Then $P^*=P$ and $Q^*:Y\to X\colon y\mapsto y$.
Finally, set $D:=\clconv \dom A.$
Then the following hold.
\begin{enumerate}
\item \label{t:georgei}
The composition $PA\colon Y\To Y$ is monotone.
\item \label{t:georgeii}
If $B\colon Y\To Y$ is a maximal monotone extension of $PA$
 (e.g., $B = \widetilde{PA}$), then
\begin{equation}
\widehat{B} := Q^*BQ+N_{D} \colon X \To X
\end{equation}
is a maximal monotone extension of $A$ and $\cldom
{\widehat{B}}= D$.
\end{enumerate}
\end{theorem}
\begin{proof}
\ref{t:georgei}:
Take $(y_1,y_1^*)$ and $(y_2,y_2^*)$ in $\gra A$.
Since $P^*=P$, $\dom A \subseteq Y$, and $A$ is monotone,
it follows that
$\scal{y_1-y_2}{Py_1^*-Py_2^*} =
\scal{Py_1-Py_2}{y_1^*-y_2^*} = \scal{y_1-y_2}{y_1^*-y_2^*}
\geq 0$.

\ref{t:georgeii}:
Let $B\colon Y\To Y$ be a maximal monotone extension of $PA$.
Using \cite[Theorem~18.8]{S98} and since
$\ospan \dom B=\ospan\dom A=Y$, we see that
\begin{equation}
\inte_{Y}{\dom B}=\inte_{Y}{\conv \dom B} = \inte_{Y}\clconv\dom B
\end{equation}
and that $\inte_Y\conv \dom A=\inte_{Y}\clconv \dom A$.  Thus
\begin{equation}
\inte_{Y}{\dom B}=\inte_{Y}{\clconv \dom B} \supseteq
\inte_{Y}{\clconv \dom A}
=\inte_{Y}{\conv\dom A}\ni 0,
\end{equation}
and hence
$0\in \inte_{Y}{\dom B}=\reli \dom B$.
Since $0\in Y\cap \reli\dom B=\ran Q\cap \reli\dom B$,
\cite[Theorem~12.43]{RW98} implies that $Q^*BQ\colon X
\To X$ is maximal monotone.
Because $\dom Q^*BQ=Q^{-1}\dom B,$
$Q0=0\in\inte_{Y}\dom B$, and $Q:X\to Y$ is continuous,
we see that $0\in \inte \dom Q^*BQ$.
On the other hand,
$0\in\reli D = \reli\dom N_D$.
Altogether,
\begin{equation}
0\in \inte \dom Q^*BQ \cap \reli\dom N_{D}.
\end{equation}
Thus, by \cite[Corollary~12.44]{RW98},
\begin{equation}
\widehat{B} = Q^*BQ + N_D \quad\text{is maximal monotone.}
\end{equation}
We shall now show that $\widehat{B}$ is an extension of $A$.
To this end, take $(a,a^*)\in\gra A$.
Since $a\in\dom A\subseteq Y$, we have $Qa = a$.
Recalling that $B$ extends $PA$, we deduce that
$Pa^* \in PAa \subseteq Ba = BQa = Q^*BQa$ and so
\begin{equation}\label{back1}
Pa^*\in Q^*BQa.
\end{equation}
On the other hand,
$a\in\dom A\subseteq D \subseteq Y$, which implies
$Y^\bot \subseteq N_D(a)$ and further
\begin{equation}\label{back2}
a^*-Pa^*\in N_{D}(a).
\end{equation}
Adding \eqref{back1} and \eqref{back2} yields
\begin{equation}
a^*\in \big(Q^*BQ+N_{D}\big)(a) = \widehat{B}a.
\end{equation}
Hence $(a,a^*)\in\gra \widehat{B}$ and we conclude that
\begin{equation} \label{e:laktheman}
\text{$\widehat{B}$ is a maximal monotone extension of $A$.}
\end{equation}
By \cite[Theorem~18.6]{S98},
$\clconv\dom \widehat{B}=\cldom\widehat{B}\subseteq \cldom
N_D= D$.
Furthermore, \eqref{e:laktheman} implies that
$\dom A\subseteq \dom\widehat{B}$ and thus
$D = \clconv\dom A \subseteq\clconv\dom\widehat{B}$.
Therefore,
$D = \cldom\widehat{B}$.
\end{proof}

Let us now illustrate and compare
Theorem~\ref{t:thursday} and Theorem~\ref{t:george}.

\begin{example}
Suppose that $X=\RR$, let $A = \Id_{\left]-1,1\right[}$
be the identity operator restricted to the open interval
$\left]-1,1\right[$,
and set $D = \clconv\dom A = [-1,1]$.
Since $\dom A$ is open, we have $A_D := A+N_D = A$.
Hence
\begin{equation}
\widetilde{A_D} = \widetilde{A} = \Id
\end{equation}
by \cite[Example~5.10]{BW07b}.
Let $B$ be an arbitrary maximal monotone
extension of $A$. Since $\gra B$ is closed,
it follows that $B$ extends $\Id|_D$.
Hence $B+N_D$ extends $\Id+N_D$, but the latter is already
maximal monotone. Thus, the operator $\widehat{B}$
provided by Theorem~\ref{t:george}.\ref{t:georgeii}
is
\begin{equation}
\widehat{B} = \Id+N_D.
\end{equation}
Therefore, Theorem~\ref{t:thursday} and Theorem~\ref{t:george}
may in general produce different maximal monotone extensions.
\end{example}

\section{Main Results}\label{main}

We now turn to the extension of (firmly) nonexpansive
operators. For recent results in this direction, all based
on Zorn's Lemma,
see \cite{RS05}, \cite{B07}, and the references therein.

\begin{theorem}[constructive extension] \label{t:1}
Let $S\subseteq X$ and let $T\colon S\to X$ be nonexpansive.
Proceed as follows.
\begin{enumerate}
\item[Step 1.] Set $F := \tfrac{1}{2}\Id + \tfrac{1}{2}T$.
\item[Step 2.] Set $A := F^{-1}-\Id$.
\item[Step 3.] Compute $\widetilde{A}$ as in Fact~\ref{f:BW}.
\item[Step 4.] Set $\widetilde{F} := (\Id + \widetilde{A})^{-1}$.
\item[Step 5.] Set $\widetilde{T} := 2\widetilde{F}-\Id$.
\end{enumerate}
Then $\widetilde{T}\colon X\to X$
is a nonexpansive extension of $T$.
\end{theorem}

\begin{proof}
Since $T$ is nonexpansive, $F$ of Step~1 is firmly
nonexpansive (Fact~\ref{f:nonexpvsfirm}).
By Fact~\ref{f:EB}.(i), $A$ is monotone.
Now Fact~\ref{f:BW} implies that $\widetilde{A}$
is a maximal monotone extension of $A$.
Hence, by Fact~\ref{f:EB}, $\widetilde{F}
\colon X\to X$ is firmly nonexpansive.
Finally, Fact~\ref{f:nonexpvsfirm} yields that $\widetilde{T}\colon X\to X$
is a nonexpansive extension of $T$.
\end{proof}

\begin{example}
Let $S\subseteq X$ and set $T=0|_S$.
Then Theorem~\ref{t:1} returns $\widetilde{T}=0$.
\end{example}
\begin{proof}
We use the notation of Theorem~\ref{t:1}.
Then $F=\tfrac{1}{2}\Id|_S$ and $A = \Id|_S$.
By \cite[Example~5.10]{BW07b}, $\widetilde{A}=\Id$.
Therefore, $\widetilde{F} = \tfrac{1}{2}\Id$ and $\widetilde{T}=0$.
\end{proof}

\begin{example} \label{ex:wednesday}
Let $d$ and $d^*$ be in $X$, and set $S=\{d\}$.
Consider the nonexpansive operator $T\colon S\to X$
given by $Td = d^*$. Then Theorem~\ref{t:1} returns
$\widetilde{T}\equiv d^*$.
\end{example}
\begin{proof}
Adopting the notation of Theorem~\ref{t:1},
we see that $\gra F = \{(d,e)\}$,
where $e=\tfrac{1}{2}(d+d^*)$.
Hence $\gra A = \{(e,d-e)\}$ and Example~\ref{ex:simplest}
yields $\widetilde{A}\colon x \mapsto x+(d-e)-e = x-d^*$.
It follows that $\widetilde{F}\colon x^* \mapsto \tfrac{1}{2}(x^* + d^*)$.
Therefore, $\widetilde{T} \equiv d^*$.
\end{proof}

We now translate the results on maximal monotone extensions
of the previous section to the setting of (firmly)
nonexpansive mappings.
The next two results are \emph{constructive} counterparts
to \cite[Corollary~5]{B07} provided some assumption on the range
is satisfied.

\begin{theorem}
Let $S\subseteq X$, and let $F\colon S\to X$ be firmly nonexpansive.
Set $D := \clconv\ran F$,
$A := F^{-1}-\Id$, and $A_D := A + N_D$.
Compute $\widetilde{A_D}$, and set
$\widetilde{F} := (\Id+\widetilde{A_D})^{-1}$.
Then the following hold.
\begin{enumerate}
\item $\widetilde{F}\colon X\to X$ is a firmly nonexpansive
extension of $F$ such that
$\clran \widetilde{F} =\clconv\ran \widetilde{F} \subseteq
\bigcap_{r \in\ran F} (r+T_D(r))$.
\item If $\conv\ran F$ is closed, then $\clran\widetilde{F} = D$.
\end{enumerate}
\end{theorem}
\begin{proof}
Since $\ran F =  \dom A$, the result is a direct consequence of
Theorem~\ref{t:thursday}.
\end{proof}

\begin{theorem}
Suppose that $X$ is finite-dimensional,
let $S\subseteq X$, and let $F\colon S\to X$ be firmly
nonexpansive.
Set $D := \clconv\ran F$ and $Y:=\ospan D$, and assume that
$0\in\reli D$.
Denote the orthogonal projector from $X$ onto $Y$ by $P$.
Let $B\colon Y\To Y$ be an arbitrary maximal monotone
extension of $P(F^{-1}-\Id)\colon Y\To Y$,
set $\widehat{B} := BP + N_D\colon X\To X$, and
$\widehat{F} := (\Id + \widehat{B})^{-1}$.
Then $\widehat{F}$ is a firmly nonexpansive extension of $F$ and
$\clran\widehat{F} = D$.
\end{theorem}
\begin{proof}
This is a direct consequence of Theorem~\ref{t:george}.
\end{proof}

We conclude this paper with a constructive version
of the original Kirszbraun-Valentine result.

\begin{theorem}[constructive Kirszbraun-Valentine extension]
\label{t:2}
Let $S\subseteq X$ and let $T\colon S\to X$
be nonexpansive.
Compute $\widetilde{T}$ as in Theorem~\ref{t:1}.
Set $D := \clconv\ran T$, and denote the orthogonal projector
onto $D$ by $P$.
Then $P\widetilde{T}\colon X\to X$
is a nonexpansive extension of $T$
such that $\clran(P\widetilde{T}) \subseteq D$.
\end{theorem}
\begin{proof}
From Theorem~\ref{t:1}, we know that
$\widetilde{T}\colon X \to X$ is a nonexpansive extension of $T$.
Since $P$ is nonexpansive, the result follows.
\end{proof}

\begin{remark}
We point out that if the given (firmly) nonexpansive
mapping is described by its graph,
then Remark~\ref{r:last1} and Remark~\ref{r:last2}
may be used to go back and forth between
the (firmly) nonexpansive mapping and its
full-domain extension, and the
corresponding monotone operators, respectively.
\end{remark}

\begin{remark}
Let $S\subseteq X$ and let $F\colon S\to X$ be firmly nonexpansive.
Then Corollary~\ref{c:neato} and
\cite[Corollary~5]{B07} --- which was proved
\emph{nonconstructively} --- guarantees the existence of
a firmly nonexpansive extension $G\colon  X\to X$ of $F$
such that $\clran G = \clconv\ran F$.
We do not know whether it is possible to use
Theorem~\ref{t:2} --- or any other constructive result ---
to obtain such a mapping;
see also \cite[Remark~8]{B07}.
\end{remark}


\section*{Acknowledgment}

Heinz Bauschke was partially supported by the Natural Sciences and
Engineering Research Council of Canada and by the Canada Research Chair
Program.
Xianfu Wang was partially
supported by the Natural Sciences and Engineering Research Council
of Canada.
\small

\bibliographystyle{amsalpha}

\end{document}